\documentclass[reqno,11pt]{amsart}

\usepackage{amssymb}
\newtheorem{Proposition}{Proposition}[section]

\newtheorem{Lemma}[Proposition]{Lemma}
\newtheorem{Theorem}[Proposition]{Theorem}
\newtheorem{Corollary}[Proposition]{Corollary}
\newcommand \starhat{\, \hat * \,}
\newcommand \la{\langle}
\newcommand \ra{\rangle}
\newcommand \valg{\Val^G}
\newcommand \valev{\Val^{+}}
\newcommand \valgst{\Val^{G*}}
\newcommand \valsm{\Val^{sm}}
\newcommand \ot{\otimes}
\newcommand\Hom{\mathbf{Hom}}
\newcommand \lcur{[\![} 
\newcommand\rcur{]\!]}
\newcommand{\ksm}{\K^{sm}}  
\newcommand{\R}{\mathbb{R}}
\newcommand{\K}{\mathcal{K}}
\newcommand{\D}{\mathbb{D}}
\newcommand\flag[2]{\left[\begin{array}{c} #1\\ #2 \end{array}
  \right]} 

\DeclareMathOperator{\kl}{Kl}
\DeclareMathOperator{\Val}{Val}
\DeclareMathOperator{\vol}{vol}
\DeclareMathOperator{\Gr}{Gr}
\DeclareMathOperator{\Dens}{Dens}

\title{Convolution of convex valuations}
\author{Andreas Bernig}
\author{Joseph H. G. Fu}
 
\email{andreas.bernig@unifr.ch}
\email{fu@math.uga.edu}
\date{\today}
\address{D\'epartement de Math\'ematiques, Chemin du Mus\'ee 23, 1700 Fribourg, Switzerland}
\address{ Department of Mathematics, 
University of Georgia, 
Athens, GA 30602, USA}

\begin{document}

\begin{abstract}
We show that the natural ``convolution" on the space of smooth, even, translation-invariant convex valuations on a euclidean space $V$, obtained by intertwining the product and the duality transform of S. Alesker (\cite{ale03b}, \cite{ale04}), may be expressed in terms of Minkowski sum. Furthermore the resulting product extends naturally to odd valuations as well. Based on this technical result we give an application to integral geometry, generalizing Hadwiger's additive kinematic formula for $SO(V)$ (\cite{schnei-wie}) to general compact groups $G \subset O(V)$ acting transitively on the sphere: it turns out that these formulas are in a natural sense dual to the usual (intersection) kinematic formulas. 
\end{abstract}

\thanks{ 
\\ Supported
  by the Schweizerischer Nationalfonds grant SNF 200020-105010/1.}
\maketitle 
\section{Introduction}

\subsection{General introduction}\label{sec:general introduction} In the series of papers \cite{ale01,ale03b,ale04a,ale04}, S. Alesker has introduced an array of fundamental operations on convex valuations,  illuminating and extending many classical constructions.  In rough terms, the space $\Val(V)$ of continuous translation-invariant valuations on a finite-dimensional euclidean space $V$ turns out to be (among other features) a commutative graded algebra with linear involution $\D$ (cf. Section \ref{sec:technical intro} below for precise definitions). 
Thinking of  $\D$ as analogous to the Fourier transform, Alesker has asked for a geometric characterization of the resulting ``convolution" obtained by intertwining the product with $\D$.
Our main technical result (Thm. \ref{thm:main} below) is that  this convolution is simply an artifact of the Minkowski sum operation.

Our main application (Thm. \ref{thm:additive kf}) may be described as follows. Alesker has shown that if $G \subset O(V)$ acts transitively on the sphere of $V$ then the dimension of the subspace $\valg\subset \Val(V)$ of $G$-invariant valuations is finite. It follows easily that if $\phi_1,\dots, \phi_N$ is a basis then given any $\psi \in \valg$ there are constants $c_{ij}= c_{ij}^\psi$ such that
\begin{equation}
\int_{\overline G} \psi(A \cap \bar g B) \, d \bar g = \sum_{i,j=1}^N c_{ij}\, \phi_i(A) \,\phi_j(B)
\end{equation}
for any compact compact convex bodies $A,B \subset V$ (here $\overline G= G \ltimes V$ is the group generated by $G$ and the translation group of $V$; the classical case $G = O(V)$ or $SO(V)$ was given originally by Blaschke). 

 By the same token there is also another, less well known, array of kinematic formulas for Minkowski addition instead of intersection. Hadwiger worked these out for $G=O(V)$ or $ SO(V)$ (cf.  \cite{schnei-wie}).
Thm. \ref{thm:additive kf} states that for general $G$ the addition and intersection formulas are simply transforms of one another under the involution $\D$, provided all elements of $\valg$ are even in the sense of section \ref{sec:technical intro} below.

The relations established in this paper also furnish the means for a fuller understanding of the integral geometry of the unitary group, originally studied in \cite{ale03b} and \cite{fu06}, including a completely explicit form of the principal kinematic formula for this group. This will be the subject of another paper.

\subsubsection{Acknowledgements} The authors wish to thank S. Alesker, L. Br\"ocker, M. Ludwig, J. Parsley and R. Varley for interesting and helpful conversations during the course of this project. This work was largely carried out during a visit in February 2006 of the first-named author to the University of Georgia, to which he wishes to express his warmest gratitude.

\subsection{Definitions and background}\label{sec:technical intro} \subsubsection{Valuations}

Throughout the paper, $V$ will denote a finite-dimensional euclidean vector
space. By $\mathcal{K}(V)$ we denote the
space of nonempty compact convex sets in $V$, endowed with the Hausdorff metric. The dense subspace $\mathcal{K}^{sm}(V)$
consists of compact convex sets with smooth and strictly convex
boundary.  

A (convex) valuation on $V$ is a complex-valued map $\phi$ on $\mathcal{K}(V)$ such
that 
\begin{displaymath}
\phi(K \cup L) + \phi(K \cap L)=\phi(K)+\phi(L)
\end{displaymath}
whenever $K,L, K \cup L \in \mathcal{K}(V)$.
The valuation $\phi$ is called translation invariant if
$\phi(K+x)=\phi(K)$ for all $x\in V, K \in \mathcal{K}(V)$. We will denote the space of continuous translation invariant valuations by $\Val(V)$ and use the abbreviated term {\it valuation}
to refer to them. 
However we note that Alesker \cite{ale05a,ale05b,alefu05,ale05d,ale06} has introduced a much broader notion of valuation that makes sense even on smooth manifolds. The restricted class of valuations considered here is in some sense an infinitesimal version of this broader notion.

A valuation $\varphi$ is said to have {\it parity} $\epsilon =\pm$ if $\varphi(-K)=  \epsilon\varphi(K)$ for all $K \in \mathcal K(V)$ (more colloquially we will call such valuations {\it even} and {\it odd} respectively); it is said to have degree $k$ if $\varphi (tK) = t^k \varphi(K)$ for all $K \in \mathcal K(V)$ and all $t\ge 0$. Put $\Val_k^\epsilon(V) \subset \Val(V)$ for the subspace of valuations of degree $k$ and parity $\epsilon$.  A theorem of P. McMullen \cite{mcmu} states that
\begin{equation} \Val(V) = \bigoplus_{0\le k \le n, \ \epsilon = \pm } \Val_{k}^\epsilon (V)
\end{equation}
where $\dim V = n$.  It follows that $\Val(V)$ carries a natural Banach space structure, with norm
\begin{equation}\label{eq:mu}
\parallel \varphi \parallel := \sup_{K \subset B} |\varphi (K)|
\end{equation}
where $B\subset V$ is a fixed bounded set with nonempty interior. Furthermore $\Val_0$ and $\Val_n$ are both one-dimensional, spanned by the Euler characteristic $\chi$ and the volume $\vol$ respectively.

 Note that $GL(V)$ acts continuously on $\Val(V)$ and on each $\Val_k^\epsilon(V)$ by $g\cdot \varphi(K):= \varphi(g^{-1}K), g \in GL(V), \varphi \in \Val(V), K \in \mathcal K(V)$.  If the map $g \mapsto g\phi$ is a smooth map $GL(V) \to \Val(V)$, then the valuation $\phi$ is said to be $GL(V)$-{\it smooth}. The subspace of $GL(V)$-smooth valuations, endowed with the $C^\infty$ topology obtained by identifying $\phi$ with the map $g \mapsto g\phi$, is
denoted by $\Val^{sm}(V)$.   

Given $A \in \K$ we define $\mu_A \in \Val$ by
\begin{equation}
\mu_A(K):= \vol (A+K),
\end{equation}
where $A+K:=\{a+k: a \in A,k \in K\}$ is the Minkowski sum. If $A \in \ksm$ then $\mu_A \in \valsm$, and if $A$ is antipodally symmetric then $\mu_A \in \Val^+$. McMullen conjectured, and Alesker \cite{ale01} proved, that
\begin{Theorem}\label{thm:mcmullen conjecture}
The linear span $\la \mu_A:A \in \K\ra$ is dense in $\Val$. The linear span $\la \mu_A:A \in \ksm\ra$ is dense in $\valsm$.
\end{Theorem}

\subsubsection{Alesker product and Crofton measures} \label{subsubsection_alesker_product}
The Alesker product is defined by  putting
\begin{equation}
\mu_A \cdot \mu_B(K)=\vol(\Delta( K) + A \times B), \quad K \in \mathcal{K},
\end{equation}
where $\Delta:V \to V \times V$ denotes the diagonal embedding, then extending to all pairs $\phi,\psi \in \valsm$ by continuity \cite{ale04}.
Using Fubini's theorem this may also be expressed
\begin{equation}\label{eq:fubini product}
\mu_A\cdot \varphi (K) := \int_V \varphi(K\cap (x - A)) \, dx,
\end{equation}
where $x-A:=\{x-a:a \in A\}$.
  This endows the space $\Val^{sm}(V)$ of smooth valuations with the structure of a commutative graded algebra. The Euler characteristic $\chi$ is the unit element for this
product. If $W \subset V$ is a subspace then the restriction map
\begin{equation}
r_W:\valsm(V) \to \valsm(W)
\end{equation}
is a homomorphism of algebras.

In the case of even valuations the product admits an alternative expression as follows. Given a linear subspace $P \subset V$, put $\pi_P:V\to P$ for the orthogonal projection.  If
$\phi \in \Val^{+,sm}_ k(V)$ then there exists a smooth measure $m_\phi$ on $\Gr_k(V)$ such that
\begin{equation}
\phi(K) = \int_{\Gr_k(V)} \vol_k(\pi_P(K)) \, dm_\phi(P).
\end{equation}
(If $2 \le k \le n-2$ then this measure is not unique.)
In particular, if $K$ is included in an affine subspace of dimension $<k$ then $\phi(K) = 0$.
We will call such a measure a {\bf Crofton measure} for $\phi$. If $\psi \in \Val^{+,sm}_ l(V)$ and $k+l \le n$ then $\phi \cdot \psi \in \Val^{+,sm}_ {k+l}(V)$ is the valuation with Crofton measure
\begin{equation}\label{eq:crofton}
m_{\phi\cdot \psi} = \sigma_*(\sin m_\phi \times m_\psi),
\end{equation}
where $\sigma: \Gr_{k}(V) \times
\Gr_{l}(V) \setminus \Delta \to \Gr_{k+l}(V)$ is the sum map $\sigma(P,Q) := P +Q$, $\Delta$ is the null set of pairs of planes $(P,Q)$ that fail to meet transversely, and $\sin:\Gr_k \times \Gr_l \to \R_{\ge 0}$ is the function determined by the relation
$$ \mathcal \vol_{k+l} (A +B) = \sin(E,F) \vol_k(A) \vol_l(B)$$
for convex sets $A \subset E , \, B \subset F$.
The valuation that results is independent of the choices of Crofton measures $m_\phi,m_\psi$. For more details, in the context of the more general non-translation-invariant case, the reader is referred to \cite{be05}.

\subsubsection{The Klain function and the $\D$ transform}
Let $\phi\in \Val_k^+(V)$. Then the restriction of $\phi$ to  a $k$-dimensional subspace $L$ is a 
multiple of the Lebesgue (or $k$-dimensional Hausdorff) measure on $L$. Put $\kl_\phi(L)$ to be the
proportionality factor; the resulting continuous function $\kl_\phi:\Gr_k(V) \to \R$ is called the {\it Klain
function} of  $\phi$. 
Klain \cite{kl00} showed that the map $ \kl:\Val^+_k(V) \to C(\Gr_k(V))$ is an injection.

The unique valuation $\phi \in \Val_k^{+}(V)$ whose Klain function
is identically $1$ is called the $k$-th intrinsic volume and denoted by
$\mu_k$. These valuations are $SO(V)$-invariant and, by Hadwiger's Characterization Theorem \cite{klro97}, span the space of all continuous
$SO(V)$-invariant valuations.  

Denote by $\perp:\Gr_k \to \Gr_{n-k}$ the orthogonal complement map. The $\D$ transform of a smooth, even and translation invariant
valuation $\phi \in
\Val_k^{+,sm}(V)$ is defined as the unique valuation
$\mathbb{D}\phi \in \Val_{n-k}^{+,sm}(V)$ with Klain function 
\begin{equation}\label{eq:klain d}
\kl_{\mathbb D\phi} = \kl_\phi \circ \perp.
\end{equation}
Equivalently, if $m_\phi$ is a Crofton measure for $\phi$ then the pushforward measure 
\begin{equation}\label{eq:crofton d}
m_{\mathbb D\phi}:= \perp_*m_\phi
\end{equation}
  is a Crofton measure for $\D \phi$.

\subsection{Statement of results}\label{sect:results}
\subsubsection{The main theorem} 
Thinking of $\mathbb D$ as analogous to the Fourier transform, it is natural to define the {\it convolution} of $\varphi, \psi \in \Val^{+,sm}(V)$ by 
\begin{equation}
\varphi * \psi := \mathbb D(\mathbb D \varphi \cdot \mathbb D \psi).
\end{equation} 
Alesker has asked for a geometric characterization of the convolution. Our main theorem gives an answer in terms of the action on the valuations $\mu_A, A\in \K(V)$, and shows that it extends in a natural way to all smooth valuations, regardless of parity.

\begin{Theorem} \label{thm:main} If $A,B$ are smooth, strictly convex and antipodally symmetric then
\begin{equation}\label{convolution_convex}
\mu_A *\mu_B = \mu_{A+B} .
\end{equation}
In fact there exists a unique bilinear, continuous product $*: \valsm\ot \valsm \to \valsm$ of degree $-n$ such that \eqref{convolution_convex} holds for general $A,B \in \ksm$. $(\valsm,+,*)$ is a commutative and associative algebra whose unit element is $\vol$. 
\end{Theorem}

Another characterization is in terms of mixed volumes:
\begin{Corollary} \label{convolution_mixed_volumes}
If $k +l \ge n$ and $A_1,\dots,A_{n-k}, B_1,\dots,B_{n-l} \in\ksm(V)$  then
\begin{equation} \label{eq:mix*}
V_{A_1,\dots,A_{n-k}} * V_{B_1,\dots,B_{n-l}} = \binom{k+l}{k}^{-1} \binom{k + l}{n}V_{A_1,\dots,A_{n-k},B_1,\dots,B_{n-l}}.
\end{equation}
\end{Corollary}

{\bf Remark.} Alesker has pointed out that Thm. \ref{thm:main} may also be stated without reference to a euclidean structure on $V$. As in the original treatment \cite{ale03b}, it is formally natural to view the transform 
$\D$ as an isomorphism $\Val_k^{+,sm}(V)\to \Val_{n-k}^{+,sm} (V^*)\otimes \Dens(V)$, where 
$\Dens(V)$ is the one-dimensional space of signed Lebesgue measures on $V$. Note that to any $0\ne m \in \Dens(V)$ there is a naturally associated $m^* \in \Dens (V^*)$ such that the product measure $m \times m^*$ on $V\times V^* \simeq T^*V$ equals the canonical Liouville measure on the cotangent bundle of $V$. Thus even in the absence of a euclidean structure it makes sense to consider the valuations $\mu_A$ as elements of $\Val(V) \otimes \Dens(V^*)$: given $A,B\in \K(V)$, put
$$\mu_A(B):= m(A+B) \cdot m^* $$
for $m,m^*$ as above, which is clearly independent of choices.
 In this language Thm. \ref{thm:main} becomes

\begin{Theorem}
There exists a continuous convolution product $*$ on
$ \Val^{sm} (V )\otimes \Dens(V^* )$, uniquely characterized by the condition
that  $\mu_A *\mu_B=\mu_{A+B}$ for any $A, B \in \K^{sm}(V )$. If $\phi,\psi \in \Val^{+,sm}(V) \otimes \Dens(V^*)$ then
$$ \phi *\psi = \D(\D^{-1}\phi \cdot \D^{-1}\psi).$$
\end{Theorem}

However, in the present paper we will continue to use the language of Thm. \ref{thm:main}.

\subsubsection{Applications} Alesker has shown that if $G\subset O(V)$ is a compact subgroup acting transitively on the sphere of $V$, then the vector space $\valg$ of  continuous, translation-invariant, $G$-invariant valuations is finite dimensional. From this it is easy to show (cf. \cite{fu06})
\begin{Proposition}\label{thm:g kinematic}
Let $\varphi_1,\dots,\varphi_N$ be a basis for $\valg$ and let $\psi \in \valg$. Then there exist constants $c^\psi_{ij}, \, 1\le i,j\le N$ such that whenever $A,B\in \K(\R^n)$
\begin{equation}\label{eq:g kinematic}
\int_{\overline G} \psi(A\cap \bar g B) \, d\bar g = \sum_{i,j} c^\psi_{ij} \varphi_i(A) \varphi_j(B).
\end{equation}
\end{Proposition}
As before, $\overline G:= G\ltimes V$ and $d\bar g$ is the Haar measure. 
The same argument from finite-dimensionality applies also to yield the following  ``additive kinematic formula."
\begin{Proposition}\label{thm:additive kinematic}
There exist constants $d^\psi_{ij}$ such that
\begin{equation}\label{eq:additive kinematic}
\int_G \psi( A + gB) \, dg = \sum_{i,j} d^\psi_{ij} \varphi_i(A)\varphi_j(B)
\end{equation}
for all $A,B \in \K$.
\end{Proposition}

We encode these facts by defining the maps $k_G, a_G:\valg\to \valg\ot\valg$ 
\begin{align}\label{def:kG}
k_G(\psi) &:= \sum_{i,j} c^\psi_{ij} \varphi_i \ot \varphi_j, \\
a_G(\psi) &:= \sum_{i,j} d^\psi_{ij} \varphi_i \ot \varphi_j,
\end{align}
where the $c^\psi_{ij}, d^\psi_{ij} $ are the constants from \eqref{eq:g kinematic} and \eqref{eq:additive kinematic}.  
 Our first application of Thm. \ref{thm:main} relates these two formulas via $\D$:
\begin{Theorem} \label{thm:additive kf} Let $G$ be as above, and suppose that $\valg \subset \valev$. Then 
\begin{equation}\label{eq:additive kf}
 a_G = (\D \ot\D)\circ k_G \circ \D.
\end{equation}
\end{Theorem} 
In other words the structure constants for the two comultiplications $k_G,a_G$ are identical provided the bases are changed appropriately. In the case $G= SO(n)$ the structure constants had been computed by Hadwiger (cf. \eqref{eq:additive so} below and \cite{schnei-wie}). 

As a second application we recall the definitions of the  ``Lefschetz operators" used in the two different forms of Alesker's Hard Lefschetz Theorem for $\Val^{+,sm}$:
\begin{align*}
\Lambda: \Val_* \to \Val_{*-1},\\
L: \Val_* \to \Val_{*+1}
\end{align*}
given by
\begin{align*}
\Lambda(\phi) (A) &:= \left.\frac{d}{dt}\right|_{t=0+} \phi(A + tB),\\
L(\phi) &:= \mu_1 \cdot \phi,
\end{align*}
where $B$ is the unit ball of $V$.

\begin{Corollary}\label{thm:Lambda/L} For all $\phi \in \valsm$, 
$$\Lambda (\phi) = 2\mu_{n-1} * \phi$$
\end{Corollary}

\begin{Corollary}
$$
\Lambda|_{\Val^{+,sm}} = 2 \D\circ L \circ \D.
$$
\end{Corollary}
It follows that the two forms of the ``hard Lefschetz theorem" for even valuations, established by Alesker in \cite{ale03b} and \cite{ale04a}, are equivalent.

\subsection{Generalities about kinematic formulas}  \label{subsec_kin_form}
We take this opportunity to clarify the relation between the product in $\valg$ and the kinematic operator $k_G$, sketched in \cite{fu06}.

Let $G$ be a subgroup of $O(V)$, acting transitively on the unit
sphere of $V$. 
We endow $G$ with the unique Haar measure of volume $1$. The
semidirect product $\bar G:=G \ltimes V$ is endowed with the product
of Haar and Lebesgue measure.  
We denote by $\Val^G(V)$ the space of $\bar G$-invariant valuations
on $V$, i.e. the space of valuations that are both translation invariant and $G$-invariant.

The best known case is $G=SO(V)$, in which case the intrinsic volumes
$\mu_0,\ldots,\mu_n$ span $\Val^{SO(V)}$. We write $\omega_n$ for the volume of
the $n$-dimensional unit ball and (following \cite{klro97}) denote by  
\begin{displaymath}
\flag{n}{k}:=\binom{n}{k} \frac{\omega_n}{\omega_k\omega_{n-k}}
\end{displaymath}
the {\it flag coefficient}. The kinematic formula for $G=SO(V)$ is (cf. \cite{klro97})

\begin{equation} \label{classical_kinematic_formula}
k_{SO(V)}(\mu_k)=\sum_{i+j=n+k} \flag{n+k}{k} \flag{n+k}{i}^{-1} \mu_i
\otimes \mu_j. 
\end{equation}

More generally, the map $k_G:\Val^G \to \Val^G \otimes \Val^G$ from \eqref{def:kG} is a cocommutative, coassociative coproduct. It is compatible with
Alesker's product in the sense that
\begin{equation}\label{multiplicativity}
k_G(\psi \cdot \phi)=(\psi
\otimes \chi) \cdot k_G(\phi)
\end{equation} 
(cf. \cite{fu06}).

Denote by 
$$
m_G: \Val^G \ot \Val^G \to \Val^G
$$
the multiplication map, and put
$$ m_G^*: \Val^{G*} \to \Val^{G*} \ot \Val^{G*}$$
for its adjoint. Recall that $\Val_n = \Val_n^G$ is spanned by a choice of a Lebesgue measure $\vol$ on $V$, where $\dim V = n$.

Let $\vol^* \in \Val^{G*}$ be the element such that
$\la \vol, \vol^* \ra =1$  and $\vol^* \perp \Val^G_k $ for $k < n$.
Recall (\cite{ale04}) that the graded algebra $\valsm$ satisfies Poincar\'e duality. Restricting to $\valg$, this means that
the map $p \in \Hom(\valg,{\valgst})$ determined by
$$ \la a,p(b)\ra\equiv\la ab,\vol^* \ra $$
is an isomorphism of vector spaces. Thus $\valg$ is structurally similar to the cohomology algebra of a compact oriented manifold $M$.  The following result states that in this picture the kinematic formula is analogous to  the linear injection $H^*(M) \to H^*(M)\otimes H^*(M)$ induced by the diagonal map and Poincar\'e duality.

\begin{Theorem}\label{thm:abstract}
\begin{equation}
\label{eq:abstract}(p\ot p )\circ k_G = m_G^* \circ p.
\end{equation}
\end{Theorem}

\begin{proof} If $a,b,c \in \valg$ then
\begin{align}\label{eq:rel1}
\la (p\ot p )\circ k_G (a),b\ot c\ra&=\la k_G(a) \cdot( b\ot c), \vol^* \ot \vol^*\ra,\\
\label{eq:rel2}
\la m_G^* \circ p(a),b\ot c\ra &= \la abc, \vol^*\ra,
\end{align}
so it is enough to show that the right hand sides of \eqref{eq:rel1} and \eqref{eq:rel2} agree.

In the case $a= \chi = 1$ this is the content of \cite{fu06}, Thm. 2.6:  if $\{d_i\}$ is any basis for $\valg$ and $\{d_i^*\}\subset{\valg}^*$ is the dual basis, then
$$ k_G(1) = \sum_i d_i \ot p^{-1}(d_i^*).$$
If this basis includes $c= d_0$,  then $\la p^{-1}(d_i^*)\cdot c,\vol^*\ra =\la d_i^*,c\ra = \delta^i_0$, so
\begin{align*}
\la k_G(1) \cdot b\ot c,\vol^*\ot \vol^*\ra &= \la(c \ot p^{-1}(c^*))\cdot (b\ot c),\vol^* \ot \vol^*\ra\\
& =\la (bc) \ot \vol,\vol^* \ot \vol^*\ra\\
&=\la bc,\vol^*\ra.
\end{align*}
Now the theorem follows from \eqref{multiplicativity}.
\end{proof}

\section{Product and convolution}

In this section we prove the main theorem \ref{thm:main}. The first part, characterizing the convolution  for even valuations in terms of the Minkowski sum, will be proved in section \ref {sect:even convolution}. The second part, extending the resulting product to all valuations regardless of parity, is proved in \ref{sect:odd convolution}.

\subsection{Convolution for even valuations}\label{sect:even convolution}
\begin{Lemma} \label{klain_function_muA}
Let $A \in \mathcal{K}^{sm}(V)$ be symmetric about the origin. Then the Klain
function of the degree $k$ component $\mu_{A,k}$ of $\mu_A$ is given by
\begin{displaymath}
\kl_{\mu_{A,k}}(L)=\mu_{n-k}(\pi_{L^\perp} A), \quad L \in \Gr_k(V). 
\end{displaymath}  
\end{Lemma}

\proof
Let $L \in \Gr_k(V)$ and let $B_L(r)$ be a ball in $L$ of radius $r$. Then
\begin{align}
\mu_A(B_L(r))&=\vol(B_L(r)+A)\\
 &= \vol (B_L(r) + \pi_{L^\perp} A)+ o(r^k)\\
& =\mu_k(B_L(r)) \mu_{n-k}(\pi_{L^\perp} A) + o(r^k)
\end{align}
as $r \to \infty$.
Since $\mu_{k}(B_L(r))=O(r^k) $, 
\begin{displaymath}
\kl_{\mu_{A,k}}(L)=\lim_{r \to \infty} \frac{\mu_A(B_L(r))}{\mu_k(B_L(r))}=\mu_{n-k}(\pi_{L^\perp}A).
\end{displaymath}   
\endproof

From this it is easy to deduce

\begin{Corollary} \label{cor:restrict}
Let $W \subset V$ be an $m$-dimensional subspace, $\mathbb{D}_W:\Val^+_* (W)\to \Val^+_{m-*} (W)$  the duality operator in
$W$,
 and $r_W:\Val(V) \to\Val(W)$ the restriction map. Then
\begin{equation}\label{eq:restriction}
r_W(\mathbb{D}\mu_{A}) =\mathbb{D}_W (\mu_{\pi_WA}).
\end{equation}
 \end{Corollary}
\endproof

Put $\cos:\Gr_k\times \Gr_k \to \R_{\ge 0}$ for the function $\cos (P,Q):= \sin(P,Q^\perp)$. Alternatively, $\vol_k(\pi_Q(A))= \cos(P,Q) \vol_k(A)$ for $A \subset P$, or vice versa. In particular $\cos(P,Q) = \cos(Q,P)$.  

\begin{Lemma} \label{lemma_D_selfadjoint}
Let $\phi,\psi \in \Val^{+,sm}_k$. 
Then 
\begin{align}
\la p(\psi), \mathbb{D} \phi\ra&= \int_{\Gr_{k}(V)} \kl_{\psi}(P) \,dm_\phi(P)\nonumber \\
\label{eq:pairing}&= \int\int \cos(P,Q) \, dm_\psi(Q) \, dm_\phi(P).
\end{align}
In particular
\begin{equation}\label{eq:symmetry}
\psi \cdot \mathbb{D} \phi=\mathbb{D} \psi \cdot \phi.
\end{equation}
\end{Lemma}

\proof This is immediate from \eqref{eq:crofton} and \eqref{eq:crofton d}.
\endproof 

\begin{Corollary} \label{cor:dxd} $$(\D\ot \D)(k_G(\chi) )= k_G(\chi). $$
\end{Corollary}
\begin{proof} Observing that $\D^*\circ p = p\circ \D$, Thm. \ref{thm:abstract} implies that
\begin{align}
(\D\ot \D)(k_G(\chi)) &= (\D \ot \D)\circ(p^{-1} \ot p^{-1}) \circ m_G^*(\vol^*)\nonumber\\
\label{eq:dxd}&= (p ^{-1}\ot p^{-1})\circ (\D ^*\ot \D^*) \circ m_G^*(\vol^*).
\end{align}
Meanwhile, given any $\phi,\psi \in \valg(V)$,
\begin{align*}
\la  \phi \ot \psi,(\D ^*\ot \D^*) \circ m_G^*(\vol^*) \ra &=   \la m_G\circ(\D \ot \D) (\phi \ot \psi),\vol^*\ra\\
&= \la \D\phi \cdot \D\psi, \vol^* \ra\\
&= \la \phi \cdot \psi, \vol^* \ra \quad \text{ (by \eqref{eq:symmetry})}\\
&= \la \phi \ot \psi, m_G^* \vol^*\ra.
\end{align*}
Thus $(\D^* \ot \D^*) \circ m_G^*(\vol^*) = m_G^*(\vol^*)$, so \eqref{eq:dxd} becomes
$$ (\D\ot \D)(k_G(\chi)) = (p ^{-1}\ot p^{-1})\circ  m_G^*(\vol^*)= k_G(\chi).$$
\end{proof}

\proof[Proof of the first part of Thm. \ref{thm:main}]

Let $A,B \in \mathcal{K}^{sm}(V)$ be centrally symmetric. Let $\mu_A=\sum_{k=0}^n \mu_{A,k}$ and $\mu_B=\sum_{k=0}^n \mu_{B,k}$ be the decompositions of $\mu_A$ and $\mu_B$à by degree of homogeneity, and let $m_{B,k}$ be a smooth Crofton measure for $\mu_{B ,k}$.  Then

\begin{align*}
\mathbb{D}\mu_{A,n-k} \cdot \mathbb{D} \mu_{B,k} & =\int_{\Gr_{k}(V)}
\kl_{\mathbb{D}\mu_{A,n-k}}(L) \, dm_{B,k}(L) \vol  \text{ by } \eqref{eq:pairing}\\ 
& =  \int_{\Gr_{k}(V)}
\vol_{k}(\pi_L(A)) \, dm_{B,k}(L) \vol \text{ by Lemma \ref{klain_function_muA}} \  \\
& =  \mu_{B,k}(A ) \ \vol .
\end{align*}
Summing over $k=0,\ldots,n$,
\begin{equation}
(\mathbb{D}\mu_A \cdot \mathbb{D} \mu_B)_n=\sum_{k=0}^n\mu_{B,k}(A) \vol=\mu_B(A)\vol=\vol(A+B) \vol=(\mathbb{D}\mu_{A+B})_n. \label{duality_highest_degree}
\end{equation}  

Now let $W \subset V$ be an $m$-dimensional subspace. Corollary \ref{cor:restrict} implies that
\begin{align*}
r_W(\mathbb{D}\mu_A \cdot \mathbb{D} \mu_B)& = r_W(\mathbb{D}\mu_A)
\cdot r_W(\mathbb{D}\mu_A)\\
& = \mathbb{D}_W \mu_{\pi_WA} \cdot \mathbb{D}_W \mu_{\pi_WB}
\end{align*}
Applying \eqref{duality_highest_degree} for the subspace $W$ and using \eqref{eq:restriction}, we conclude that the degree $m$ component of this valuation is
\begin{equation}
(r_W(\mathbb{D}\mu_A \cdot \mathbb{D} \mu_B))_m=(\mathbb{D}_W \mu_{\pi_W(A+B)})_m\\
 = (r_W(\mathbb{D}\mu_{A+B}))_m.
\end{equation}
Since this holds for all $W \in \Gr_m(V)$, we deduce
that the Klain functions of $(\mathbb{D}\mu_A \cdot \mathbb{D}
\mu_B)_m$ and $(\mathbb{D}\mu_{A+B})_m$ coincide for all
$m=0,\ldots,n$. Therefore $\mathbb{D}\mu_A \cdot \mathbb{D}
\mu_B=\mathbb{D}\mu_{A+B}$.
\endproof

\subsection{Convolution for odd valuations}\label{sect:odd convolution}
As we have mentioned above, the convolution extends in a natural way to include also the (smooth) odd valuations. Let us fix an orientation of $V$.  
 
\subsubsection{Representation of valuations by means of currents} 
Recall that a differential form on the sphere bundle $SV$ is said to be {\it vertical} if it annihilates the contact distribution of $SV$. For $0 \le k \le n-1$, put $\Omega^{V}_k$ for the space of smooth, translation-invariant differential forms $\beta$ of bidegree $(k,n-k-1)$ on $SV$ such that $d\beta$ is vertical. Put $\Omega^V:= \bigoplus_{k=0}^{n-1} \Omega^{V}_k$. 
 
The next statement is a special case of Theorem 5.2.1.~ of \cite{ale05a} and of the main theorem of \cite{bebr05}. Let $N(K) \in \mathbb{I}_{n-1}(SV)$ be the normal cycle of $K$ \cite{zaeh86,alefu05}. 
\begin{Lemma} \label{lm:kernel} If $k \ne n$ then the map $\nu:\Omega^{V}_k\to \Val^{sm}_k(V)$, given by
\begin{equation}\label{eq:n rep}
\nu(\beta)(K) := \int_{N(K)} \beta,
\end{equation}
is surjective, with kernel
\begin{equation}
\ker \nu = \{\beta\in \Omega^{V}: \beta\text{ is  exact}\}.
\end{equation}
\end{Lemma}

{\bf Definition.} Let $\pi_1:SV \to V$ and $\pi_2:SV \to S(V)$ denote the canonical projections, where $S(V)$ is the unit sphere of $V$. Let $*_V$ be the Hodge star on the space $\Omega^*(V)$ of differential forms on $V$, and let $*_1$ be the
linear operator on $\Omega^*(SV)$ which is uniquely defined by  
\begin{align}\label{eq:sign convention} 
*_1(\pi_1^* \gamma_1 \wedge \pi_2^* \gamma_2)&=(-1)^{\binom{n-\deg \gamma_1}{2}} \pi_1^* (*_V \gamma_1) \wedge
 \pi_2^* \gamma_2, \\
 \nonumber  \gamma_1 \in \Omega^*(V),&\quad\gamma_2 \in \Omega^*(S(V)).
\end{align} 
We define $*_1$ on $\Omega^*(TV)$ in a similar way. 
 
Note that $d *_1=(-1)^n *_1 d$ on translation invariant forms and that $*_1$ is, up to a sign, an involution. 

If $\beta\in \Omega^{V}_k, \gamma \in \Omega^{V}_l$ put
\begin{align}
\beta\starhat  \gamma &:= (2n-k-l)^{-1}*_1^{-1} \left((n-k) *_1\beta \wedge *_1d\gamma + (n-l) *_1 \gamma \wedge *_1 d\beta\right) \nonumber\\
\label{eq:starhat}&\equiv *_1^{-1} (*_1 \beta \wedge *_1 d\gamma) \mod \ker \nu. 
\end{align}

\begin{Proposition} Extending by bilinearity, $\starhat$  is a continuous, commutative, associative product on $\Omega^V$ of degree $-n$, which descends to a continuous, commutative, associative product $*$ on $\valsm$ of degree $-n$ by taking
\begin{align}
\label{eq:val *} \nu(\beta) *\nu(\gamma)&:= \nu(\beta \starhat \gamma), \\
\label{eq:vol *} \vol * \, \phi &:= \phi.
\end{align}
\end{Proposition}

The reader will observe that the definition in \eqref{eq:val *}, \eqref{eq:vol *} is a clear abuse of notation, since for even valuations we defined $*$ differently in the last section. However, we will see shortly (Prop. \ref{prop:*=*} below) that the two definitions agree in this case, so we hope the reader will tolerate this momentary formal ambiguity. Note also that $*$, unlike $\starhat$, does not depend on the orientation of $V$. 

\proof Commutativity and associativity may be verified via straightforward computations, taking into account that $\deg *_1 \beta$ is always odd for $\beta \in \Omega^V$.

To see that $*$ in \eqref{eq:val *} is well defined it is enough to show that if $\beta \in \Omega^V$ is exact then so is $\beta \starhat \gamma$ for all $\gamma \in\Omega^V$. However this follows at once from \eqref{eq:starhat}.

To prove continuity, let $\phi_1^j \to \phi_1, \phi_2^j \to \phi_2$ in $\valsm$. We may assume that the degrees of all of these valuations are $<n$. Since the map $\Omega^V \to \valsm_{<n} $ is surjective, and is obviously continuous with respect to the $C^\infty$ topology on $\Omega^V$, the open mapping theorem implies that we can choose sequences
$\beta_1^j, \beta_2^j\in \Omega^V$ representing $\phi_1^j, \phi_2^j$ and converging to
$\beta_1, \beta_2$ in the $C^\infty$ topology. Then $\beta_1^j \starhat
\beta_2^j \to \beta_1 \starhat \beta_2$ in $C^\infty$, and thus $\phi_1^j * \phi_2^j \to
\phi_1 * \phi_2$.
\endproof

It remains to show

\begin{Proposition}\label{prop:*=*} If $*$ is defined as in \eqref{eq:val *}, \eqref{eq:vol *} then
\begin{equation}\label{eq:starhat minkowski} \mu_A *\mu_B = \mu_{A+B}
\end{equation}
for all $A,B \in\ksm$.
\end{Proposition}
\proof It will be convenient to use a variation on the representation \eqref{eq:n rep} of valuations as integrals over the normal cycle $N$. Thinking of $N$ as analogous to the manifold of unit normal vectors of a submanifold, we use instead the analogue of the bundle of unit balls:
\begin{equation}
N_1(K):= \vec N(K) \, \llcorner\, (V \times B(0,1)), \ K \in \K,
\end{equation}
where $\vec N$ is the image under the identification $TV \simeq T^*V$, induced by the euclidean metric, of the conormal cycle $\vec N^*$ of \cite{alefu05}.

 Let us put $r(x,y) := |y|$ and $p:TV \setminus V \to SV, (x,y) \mapsto \left(x,\frac{y}{r}\right)$. Given $\beta \in \Omega^{V}_k$, let 
$$ \tilde \beta := \begin{cases} d\left( r^{n-k}\,p^*\beta\right), &\text{ on } TV \setminus V\\
	0 &\text{ on } V.
\end{cases}
$$ 
Since $\tilde \beta$ is the differential of a Lipschitz form it defines a flat cochain in the sense of \cite{wh57} (cf. also \cite{fe69}, 4.1.19), and the usual form of Stokes' theorem applies to give
\begin{equation}\label{eq:n1 rep}
\int_{N_1(K)} \tilde \beta = \int_{\partial N_1(K)} r^{n-k}p^*\beta = \int_{N(K)} \beta = \nu(\beta)(K)
\end{equation}
for all $K \in \K$. Furthermore the volume valuation is given by
$$
\vol(K) = \int_{N_1(K)} \pi_1^*(d\vol_V).
$$

Put $ \tilde \Omega^{V}_k: = \{\tilde \beta: \beta \in \Omega^{V}_k\}$ for $k\ne n$, and $\tilde \Omega_{n}^{V}:= \la \pi_1^* \vol_V\ra$; put $\tilde \Omega^V: =\bigoplus_{k=0}^{n}\tilde \Omega_k^V$. Put $\tilde \nu(\tilde \beta) $ for the left hand side of \eqref{eq:n1 rep}. Thus $\tilde \nu:\tilde \Omega^V \to \valsm$ is surjective.  
For $\phi, \psi \in \tilde \Omega^V$ we put
$$\phi\, \tilde * \, \psi : =   *_1^{-1}(*_1 \phi\wedge *_1\psi).$$
Unwinding the definitions, if $\beta,\gamma \in \Omega^V$ then
\begin{equation}
\tilde \beta \,\tilde * \, \tilde \gamma = \widetilde {\beta \starhat \gamma}.
\end{equation}
In particular, by \eqref{eq:val *}, \eqref{eq:vol *} 
\begin{equation}\label{eq:tilde product}
\tilde \nu(\phi) * \tilde \nu(\psi) = \tilde \nu(\phi \, \tilde * \, \psi)
\end{equation}
for all $\phi, \psi \in\tilde \Omega ^V$.

\begin{Lemma}\label{lm:lifts} Let $\theta$ be a translation-invariant form of bidegree $(k,n-k)$ on $TV$. Then $\theta \in \tilde\Omega^{V}_k$ iff $\theta$ is closed, homogeneous of degree $n-k$, and 
$\theta|_{SV}$ is vertical.
\end{Lemma}
\proof[Proof of lemma]  Suppose $k<n$ and that $\theta$ has these properties. By homogeneity,
$$ \theta = d(r^{n-k}) \wedge p^* \beta + r^{n-k}\, p^*\gamma$$
for some forms $\beta,\gamma$ on $SV$, where $\beta$ has bidegree $(k,n-k-1)$. If $d\theta = 0$ then $\gamma = d\beta$, whence $\theta = d\left(r^{n-k} p^*\beta\right)$. The restriction of $\theta$ to $SV$ is then equal to $d\beta$, which is therefore vertical, from which we conclude that $\beta \in \Omega^{V}_k$ and $\theta=\tilde \beta$.

The converse, and the case $k=n$, are trivial.
\endproof

Returning to the proof of \eqref{eq:starhat minkowski}, given $A \in \ksm$ we put $h_A:V\to \R$ for the support function $h_A(y):= \max_{x\in A} \la y,x\ra$ of $A$,  put
$$
\eta= \eta(y) := \begin{cases}r \nabla h_A(y), & y\ne 0,\\
0, & y = 0
\end{cases}
$$ 
and define
the Lipschitz map $G_A: TV \to V$ by
$$ 
 G_A(x,y) := x + \eta(y).
$$
We claim that for any $K \in \K$
\begin{equation}\label{eq:ga}
G_{A*} N_1(K) = \lcur K + A\rcur.
\end{equation}
Since $G_{A*}:\mathbb I_*(TV)\to \mathbb I_*(V) $ and $N_1: \K \to \mathbb I_{n}(TV)$ are continuous (cf. \cite{fe69}, \cite{alefu05}), it is enough to prove this in the case that $K \in \ksm$. Furthermore, since the left-hand side represents an integral current of top dimension with compact support,  it is enough to show that $G_{A*}N(K) = \lcur\partial( K+ A)\rcur.$ In this case $N(K)$ is given by integration over $\{(x,n_x): x\in \partial K\}$, where $n_x$ is the unit outward normal to $K$ at $x$. Thus $G_{A*}N(K)$ is given by integration over the smooth hypersurface $H:= \{x + \nabla h_A(n_x):x\in \partial K\}$. It is easy to verify that $n_x$ is normal to $H$ at $G_A(x,n_x)$, and that
\begin{align*}
\la n_x, G_A(x,n_x)\ra &= \la n_x, x\ra + \la n_x, \nabla h_A(n_x)\ra\\
&= h_K(n_x) + h_A(n_x) \\
&= h_{K+A}(n_x).
\end{align*}
This is sufficient to characterize $H$ as $\partial (K + A)$, as claimed.

Therefore
\begin{equation} \label{codisc_bundle_mink}
\vol(K+A)=\int_{G_{A*}N_1(K)}d\vol_V=\int_{N_1(K)}
G_A^* (d\vol_V).
\end{equation}
We claim that  $\theta_A:= G_A^*(d\vol_V) \in \tilde \Omega^V$. Appealing to Lemma \ref{lm:lifts}, it is clear that in the decomposition of $\theta_A$ by bidegree the components are all closed and of the correct homogeneity. It remains only to show that $\theta_A|_{SV}$ is vertical, i.e. that $\alpha \wedge \theta_A|_{SV} = 0$, where $\alpha(x,y) = \sum_{i=1}^n y_i \,dx_i$ is the canonical 1-form.

Since $h_A$ is homogeneous of degree 1 it follows that $\sum_i y_i d\eta_i |_{SV}=0$, and we establish the claim by computing
 \begin{align*}
\alpha \wedge \theta_A|_{SV} &=\left (\sum_i y_i dx_i\right) \wedge G_A^* (d\vol_V)|_{SV} \\
&= \left (\sum_i y_i \left(dx_i+ d\eta_i\right) \right) \wedge G_A^* (d\vol_V)|_{SV} \\
&= \left (\sum_i y_i (G_A^* dx_i)\right) \wedge G_A^* (d\vol_V)|_{SV} \\
&= \sum_i y_i G_A^* \left(dx_i\wedge d\vol_V\right) |_{SV} \\
&=0.
\end{align*}

Thus $\mu_A= \tilde \nu (\theta_A)$, and since
\begin{displaymath}
\theta_A:=\bigwedge_{i=1}^n  G_A^* dx_i=\bigwedge_{i=1}^n (dx_i+d\eta_i),
\end{displaymath}
a straightforward computation using \eqref{eq:sign convention} reveals that 
\begin{equation}
*_1 \theta_A=\bigwedge_{i=1}^n (1+dx_i \wedge d\eta_i).
\end{equation} 
Therefore if $B \in \ksm$ and $\zeta:= r\nabla h_B(y)$ then
\begin{align*}
*_1 (\theta_A \, \tilde* \,\theta_B) & =*_1 \theta_A \wedge *_1\theta_B \\
& =\bigwedge_{i=1}^n
 (1+dx_i \wedge d\eta_i ) \wedge (1+dx_i \wedge
 d\zeta_i) \\
& = \bigwedge_{i=1}^n
 (1+dx_i \wedge (d\eta_i + d\zeta_i))\\
& = *_1 \theta_{A+B}
\end{align*}
since $\eta + \zeta = r(\nabla h_A + \nabla h_B) = r\nabla h_{A+B}$.
Therefore $\theta_A\, \tilde *\, \theta_B=\theta_{A+B}$, which with \eqref{eq:tilde product} completes the proof. 
\endproof

Let us denote by $\tilde \sigma:\valsm \to \valsm$ the involution defined by $\tilde \sigma \phi(K):=\phi(-K)$ ($\tilde \sigma$ agrees up to a sign with the Euler-Verdier involution $\sigma$ introduced in \cite{ale05b}). Note that even valuations are fixed by $\tilde \sigma$ and that $\tilde \sigma \mu_A=\mu_{-A}$. Relation \eqref{eq:symmetry} may be rewritten as $\phi \cdot \psi=\mathbb{D} (\phi * \tilde \sigma \psi)$ for $\phi,\psi$ even and of complementary degrees. In this form, this formula holds for all smooth valuations:
\begin{Corollary}
For $\phi \in \Val_k^{sm}(V), \psi \in \Val^{sm}_{n-k}(V)$, 
\begin{equation} \label{equivalence_of_inner_products}
\phi \cdot \psi=\mathbb{D} (\phi * \tilde \sigma \psi).
\end{equation}
\end{Corollary}  

\proof
An equivalent formulation of \eqref{equivalence_of_inner_products} is that for arbitrary $\phi,\psi \in \Val^{sm}(V)$, the highest degree part of $\phi \cdot \psi$ equals the lowest degree part of $\phi * \tilde \sigma \psi$. By continuity and linearity it suffices to show this for $\phi=\mu_A$ and $\psi=\mu_B$ with $A,B \in \mathcal{K}^{sm}$. 
The lowest degree component of $\mu_A * \tilde \sigma \mu_B=\mu_{A-B}$ obviously equals $\vol(A-B) \chi$. 
On the other hand, by \eqref{eq:fubini product} the highest degree component of $\mu_A \cdot \mu_B$ is  $\vol(A-B) \vol$.
  
\endproof

\subsection{An alternative formulation} 
\proof[Proof of Corollary \ref{convolution_mixed_volumes}]
By definition of the mixed volume, 

\begin{align*}
V_{A_1,\ldots,A_{n-k}} & =\frac{k!}{n!} \left. \frac{\partial^{n-k}}{\partial t_1 \ldots
    \partial t_{n-k}}\right|_{t=0} \mu_{\sum_{i=1}^{n-k} t_i A_i}, \\
V_{B_1,\ldots,B_{n-l}} & =\frac{l!}{n!} \left. \frac{\partial^{n-l}}{\partial s_1 \ldots
    \partial s_{n-l}}\right|_{s=0} \mu_{\sum_{j=1}^{n-l} s_j B_j},
\end{align*}
from which we deduce (using  \eqref{convolution_convex} and continuity) that 
\begin{align*}
V_{A_1,\ldots,A_{n-k}} * V_{B_1,\ldots,B_{n-l}} & =\frac{k!l!}{n!^2} \left. \frac{\partial^{n-k}}{\partial t_1 \ldots
    \partial t_{n-k}}\right|_{t=0} \left. \frac{\partial^{n-l}}{\partial s_1 \ldots
    \partial s_{n-l}}\right|_{s=0} \\
& \quad \mu_{\sum_{i=1}^{n-k} t_i A_i +
    \sum_{j=1}^{n-l} s_j B_j}\\
& = \binom{k+l}{k}^{-1} \binom{k+l}{n} V_{A_1,\ldots,A_{n-k},B_1,\ldots,B_{n-l}}.  
\end{align*}
\endproof


\section{Proofs of Theorem \ref{thm:additive kf} and Corollary \ref{thm:Lambda/L}}
\proof[Proof of Thm. \ref{thm:additive kf}] We claim first that
\begin{equation} \label{coassociative}
a_G(\psi * \phi)= (\psi \otimes \vol) * a_G(\phi) = (\vol \otimes \psi) * a_G(\phi),
\end{equation}
for $\psi,\phi \in \Val^G(V)$. 
The valuations $\mu_L^G$ defined by
\begin{displaymath}
\mu_L^G(K)=\int_G \vol(K+gL)\, dg,
\end{displaymath}
 $L \in \mathcal{K}^{sm}(V)$, span $\Val^G(V)$.  It thus suffices to show  \eqref{coassociative} for $\psi=\mu_L^G$, in which case we compute 
\begin{align*}
a_G(\mu_L^G * \phi)(A,B) & = \int_G \mu_L^G * \phi(A+gB) \, dg \\
& = \int_G \int_G \phi(A+gB+hL)\, dh\, dg\\
& = \int_G \int_G \phi(A+gB+hL)\, dg\, dh\\
& = \int_G \sum_{i,j=1}^N d_{i,j}^\phi \varphi_i(A+hL) \varphi_j(B)\, dh\\
& = \sum_{i,j=1}^N d_{i,j}^\phi (\mu_L^G * \varphi_i)(A) (\vol * \varphi_j)(B)\\
& = (\mu_L^G \otimes \vol) * a_G(\phi)(A,B),
\end{align*}
which is the first equation of \eqref{coassociative}, and the second follows similarly. 

To prove \eqref{eq:additive kf} we compute
\begin{multline*}
a_G(\vol)(A,B) =\int_G \vol(A+gB) dg = \int_{\bar G} \chi(A \cap gB) dg= k_G(\chi)(A,B)= \\
= (\mathbb{D} \otimes \mathbb{D}) (k_G(\chi))(A,B)
\end{multline*}
by Corollary \ref{cor:dxd}. In other words, both sides of 
 \eqref{eq:additive kf} give the same result when evaluated on the valuation $\vol$. Now
\eqref{coassociative} and \eqref{multiplicativity} yield
\begin{align*}
a_G(\phi) & = a_G(\phi * \vol)\\
& = (\phi \otimes \vol) * a_G(\vol)\\
& = (\mathbb{D} \otimes \mathbb{D})((\mathbb{D}\phi \otimes  \chi) \cdot k_G(\chi))\\
& = (\mathbb{D} \otimes \mathbb{D})(k_G(\mathbb{D} \phi)),
\end{align*}
as claimed.
\endproof

{\bf Remarks.} (1) 
Let $c_G: \valg\ot\valg \to \valg$ denote the convolution map, determined by
$$c_G (\varphi\ot\psi) := \varphi * \psi.$$ 
Using the fact that the relation \eqref{eq:symmetry} may be stated
\begin{equation}
p \circ \D = \D^* \circ p,
\end{equation}
it is easy to deduce that, analogously to Thm. \ref{thm:abstract},
\begin{equation}
(p\ot p )\circ a_G = c_G^* \circ p.
\end{equation}

(2) Thm. \ref{thm:additive kf} says that the coefficients for the additive kinematic
formulas in a given basis are the same as the coefficients of the
ordinary kinematic formulas for the Alesker dual basis. For instance, if
$G=SO(V)$ then
\begin{align*}
a_{SO(V)}(\mu_k) & =(\mathbb{D} \otimes \mathbb{D})
 k_{SO(V)}(\mu_{n-k})\\
& =\flag{2n-k}{n-k} \sum_{i+j=2n-k}
 \flag{2n-k}{i}^{-1} \mathbb{D}(\mu_i) \otimes
\mathbb{D}(\mu_j)\\
& =\flag{2n-k}{n-k} \sum_{i+j=k} \flag{2n-k}{n-i}^{-1} \mu_i \otimes \mu_j,
\end{align*} 
 yielding the following formula of Hadwiger (cf. \cite{schnei-wie}):
\begin{equation}\label{eq:additive so}
\int_{SO(V)} \mu_k(A+gB) dg = \flag{2n-k}{n-k} \sum_{i+j=k} 
\flag{2n-k}{n-i}^{-1} \mu_i(A) \mu_j(B)
\end{equation} 
for $A,B \in \mathcal{K}(V)$.

\proof[Proof of Corollary \ref{thm:Lambda/L}]
Let $B$ be the unit ball and $r>0$. By Steiner's formula,  
\begin{displaymath}
\mu_{rB}(K)=\vol(K+rB)=\sum_{i=0}^n \mu_{n-i}(K) \omega_i r^i,
\end{displaymath}
from which we deduce that $\left. \frac{d}{dr} \right|_{r=0}
\mu_{rB}=2\mu_{n-1}$. For $\phi \in \Val^{sm}(V)$ it follows that 
\begin{equation}\label{eq:d and *}
\Lambda \phi= \left. \frac{d}{dr}\right|_{r=0} \phi(\cdot + rB) = \left. \frac{d}{dr}\right|_{r=0} \mu_{rB} * \phi = 2\mu_{n-1} * \phi.
\end{equation}
\endproof


\end{document}